\documentclass{amsart}
\usepackage{epsf}

\newtheorem{prop}{Proposition:}
\newtheorem{de}{Definition:}
\newtheorem{rk}{Remark:}
\newtheorem{lem}{Lemma:}

\newtheorem{ex}{Example:}

\newcommand{\he}{Hermitian}
\newcommand{\oo}{(1,1)-geodesic}
\newcommand{\ri}{Riemannian}
\newcommand{\ph}{pluriharmonic}
\renewcommand{\iff}{if and only if}
\newcommand{\ha}{harmonic}
\newcommand{\ho}{holomorphic}
\newcommand{\ka}{K{\"a}hler}
\newcommand{\mfd}{manifold}

\newcommand{\sn}{{\mathbb S}}

\newcommand{\cn}{{\mathbb C}}

\newcommand{\dif}[2]{{\frac{\partial{#1}}{\partial{#2}}}}
\newcommand{\ssc}[1]{{\scriptscriptstyle{#1}}}

\DeclareMathOperator{\g}{grad}

\begin{document}

\begin{abstract}
Pluriharmonic maps form an important class of {\ha} maps which includes {\ho} maps. We study their morphisms, in particular the inter-relationships between {\oo}, {\ph} and $\pm${\ho} maps. Then we characterise {\ph} morphisms between {\he} {\mfd}s. We make a special study of the situation where the target is {\ka}, {\ph} morphisms being particularly well understood for this case.
\end{abstract}

\title{\bf{Pluriharmonic Morphisms}}
\author{E.~Loubeau}
\subjclass{Primary: 32F05 Secondary: 31C10, 53C55, 58E20}
\keywords{pluriharmonic maps, $(1,1)$-geodesic maps, harmonic maps}
\curraddr{Department of Pure Mathematics \\
            University of Leeds \\ 
           LS2 9JT Leeds \\
                 U.K.}
\address{(from 1st September 1996)\\ 
Universit{\'e} de Bretagne Occidentale \\
UFR Sciences et Techniques\\
Departement de Mathematiques\\
6, avenue Victor Le Gorgeu\\
BP 809\\
29285 BREST Cedex\\
France}
\email{pmt5el@amsta.leeds.ac.uk}
\maketitle

\section{Introduction}
One of the most important properties of {\ha} maps is that, when mapping from a two-dimensional domain, the energy functional is invariant under conformal changes of the metric on that domain. We can therefore talk of {\em {\ha} maps from Riemann surfaces} and their study has been prolific (cf.~\cite{EelLem78,EelLem88}). As any {\he} structure on a Riemann surface is {\ka}, {\ha} maps from a Riemann surface are {\ha} with respect to any {\he} structure on the domain. Maps which generalise this property to higher dimensions are called {\em pluri{\ha} maps}. Harmonic morphisms between {\ri} {\mfd}s are maps which pull back germs of {\ha} functions to germs of {\ha} functions. They were characterised independently by B.~Fuglede in~\cite{Fug78} and T.~Ishihara in~\cite{Ish79A} as horizontally weakly conformal {\ha} maps. A consequence of this result is that {\ha} morphisms pull back {\ha} {\em maps} to {\ha} maps as well. 
\newline Our aim is to generalise the idea of {\ha} morphism to maps which pull back {\ph} functions to {\ph} functions; we shall call such mappings {\em {\ph} morphisms}. First, we recall some types of {\he} structures, concentrating on the geometrical interpretations of their definitions.
These results enable us to compare how and when the notions of {\oo}, {\ph}ity and {\ho}ity overlap. In particular, we show that, for maps between {\he} {\mfd}s, the classes of {\ph} maps and $\pm${\ho} maps coincide precisely when the complex structure carried by the target {\mfd} is {\ka}. We find that a necessary and sufficient condition for a map between {\he} {\mfd}s to be a {\ph} morphism is that it be a $\pm${\ho} {\ph} map. We remark that such a map will pull back {\ph} maps to {\ph} maps. Combining this characterisation of {\ph} morphisms with an earlier result, we observe that {\ph} morphisms from {\he} to {\ka} {\mfd}s are exactly the $\pm${\ho} maps.
\newline I would like to thank J.~C.~Wood for painstakingly supervising this work.

\section{Almost Hermitian Manifolds}

\begin{figure}
\epsfxsize=3.5in
\epsfbox[0 0 570 405]{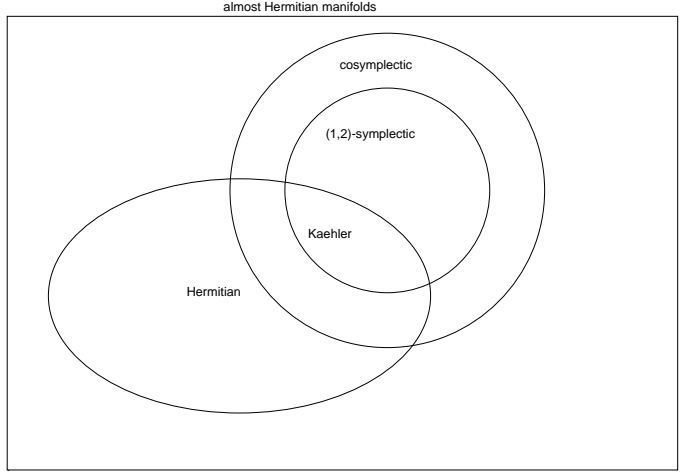}
\caption{Intersection of the structures}
\label{plot}
\end{figure}

Recall that an almost {\he} {\mfd} $(M,J,g)$ is an even dimensional manifold $M^{2n}$ equipped with an almost complex structure $J$ and a Hermitian metric $g$. On such a {\mfd} we consider an adapted unitary frame $\{e_{1},Je_{1},\dots ,e_{n},Je_{n}\}$. The complexified tangent space $T^{\ssc\cn}M$ splits into $T^{(1,0)}M$ and $T^{(0,1)}M$ which are the eigenspaces of the endomorphism $J$ with eigenvalues $i$ and $-i$, respectively.
We define $\{\theta_{j} = \frac{1}{\sqrt{2}}( e_{j} - iJe_{j} )\}_{j=1,\dots ,n}$ to be a unitary frame for $T^{(1,0)}M$ and 
$\{\theta_{\overline{\jmath}} = \frac{1}{\sqrt{2}}( e_{j} + iJe_{j} )\}_{j=1,\dots ,n}$ to be a unitary frame for $T^{(0,1)}M$.
Let $\{\Theta^{i}, \Theta^{\overline{\imath}} \}_{i=1,\dots,n}$ denotes the dual frame of $\{\theta_{i}, \theta_{\overline{\imath}}\}$, $\{\Theta^{i} \}_{i=1,\dots,n}$ (resp. $\{ \Theta^{\overline{\imath}} \}_{i=1,\dots,n}$) spans $T^{(1,0)\ast}M$ (resp. $T^{(0,1)\ast}M$).
We shall denote by $g_{i\overline{\jmath}}$ the quantity $g(\theta_{i},\theta_{\overline{\jmath}})$. 
\newline The fundamental 2-form $\omega$ is defined as:
\begin{equation}
\omega = \frac{i}{2\pi} \: g_{\lambda\overline{\mu}} \: \Theta^{\lambda}\wedge \Theta^{\overline{\mu}}. \label{defomeg}
\end{equation}
\begin{de} \label{IVde}
An almost Hermitian manifold $(M,J,g)$ is said to be:
\begin{enumerate}
   \item {  \em (1,2)-symplectic} if $(d\omega)^{1,2} =0$, i.e.  
$d\omega(\alpha,\beta,\overline{\gamma}) =0$ for all $\alpha,\beta,\gamma \in T^{(1,0)}M$.
   \item { \em {\ka}} if $\omega$ is closed, i.e. $d\omega =0$ and the almost complex structure is integrable.
   \item { \em cosymplectic} if  $d^{\ast}\omega =0$, where  $d^{\ast}$ is the formal adjoint of $d$. \label{cosymp}
\end{enumerate}
\end{de}
\begin{rk} \label{remarkka}
The names {\em quasi-{\ka}}
or {\em ${}^{\ast}0$} and {\em semi-{\ka}} or {\em almost semi-{\ka}} are sometimes used for {\em (1,2)-symplectic} and {\em cosymplectic}, respectively (cf.~\cite{Sal85}).
A list of sixteen different structures possible on an almost complex manifold as well as their classification can be found in Gray and Hervella~\cite{GrayHerv}.
\newline (1,2)-Symplectic and {\ka} manifolds can also be characterised in a more geometrical way:

If we denote by $\nabla$ the Levi-Civita connection then an almost {\he} {\mfd} $(M,J,g)$ is (1,2)-symplectic if and only if

$$ \nabla_{\ssc X} T^{(1,0)}M \subseteq T^{(1,0)}M \quad \forall X\in T^{(0,1)}M , $$
whilst it is {\ka} if and only if

$$ \nabla_{\ssc X} T^{(1,0)}M \subseteq T^{(1,0)}M \quad \forall X\in TM .$$
 
Figure~\ref{plot} on page~\pageref{plot} (based on that in~\cite{Sal85}) shows how those different structures intersect one another, in particular that an integrable (1,2)-symplectic manifold is {\ka}.
\end{rk}

The following standard proposition interprets the above definitions in terms of the metric $g$ and the Christoffel symbols $\Gamma_{\ssc{J}\ssc{K}}^{\ssc{I}}$ with respect to a unitary frame.
\begin{prop}  \cite{Lic70} \label{link}
An almost Hermitian manifold $(M,J,g)$ is
   \begin{enumerate}
         \item integrable if and only if $\Gamma_{\lambda\mu}^{\overline{\alpha}} = 0 \quad \forall \,\alpha ,\lambda ,\mu = 1,\dots,n$.
         \item  (1,2)-symplectic if and only if  
$\Gamma_{\overline{\lambda}\mu}^{\overline{\alpha}} = 0 \quad \forall \,\alpha,\lambda ,\mu$.
         \item {\ka} if and only if $J$ is integrable (equiv. $\Gamma_{\lambda\mu}^{\overline{\alpha}} = 0 \quad \forall \,\alpha ,\lambda ,\mu$), and (1,2)-symplectic (equiv. 
$\Gamma_{\overline{\mu}\lambda}^{\overline\alpha} =0 \quad \forall \,\alpha ,\lambda ,\mu $) i.e. if and only if the only possibly non-zero Christoffel symbols are $\Gamma_{\mu\lambda}^{\alpha}$ and $\Gamma_{\overline{\mu}\, \overline{\lambda}}^{\overline\alpha}$.
         \item  cosymplectic if and only if  $g^{\overline{\mu}\lambda} \Gamma^{\alpha}_{{\overline\mu}{\lambda}} =0$ for all $\alpha\in \{1,\ldots ,n\}$.                                                                                                                                                                                                                                                                                                                                                                                                                                                                                                                                                                                                                                                                       \end{enumerate}
\end{prop}

\section{Pluriharmonic Morphisms}

Let $N^{n}$ be an almost Hermitian manifold with a local adapted frame $(\vartheta_{1},\vartheta_{\overline{1}},\dots ,\vartheta_{n}, \vartheta_{\overline{n}})$ 
and $X^{r}$ a {\ri} {\mfd} with local coordinates $(\varpi_{1},\dots ,\varpi_{r})$.
We denote by $\nabla$ the Riemannian connection on $X$ and by  ${}^{\ssc N}\Gamma^{\scriptscriptstyle{C}}_{{\scriptscriptstyle{A}}{\scriptscriptstyle{B}}}$ and ${}^{\ssc X}\Gamma^{\scriptscriptstyle{C}}_{{\scriptscriptstyle{A}}{\scriptscriptstyle{B}}}$ the Christoffel symbols on $N$ and $X$

\begin{de}
A map $\phi$ from an almost {\he} {\mfd} $N^{n}$ to a {\ri} {\mfd} $X^{r}$ is called {\em{ {\oo}}} (see \cite{EelLem88}) if the (1,1)-part of the second fundamental form of $\phi$, $({}^{\phi}\nabla d\phi )^{(1,1)}$, is zero, i.e.:
\begin{equation}
\label{glob}
 {}^{\phi}\nabla d\phi (Z,{\overline{W}}) = 0 \qquad \forall Z,W\in T^{(1,0)}N ,
\end{equation}
where ${}^{\phi}\nabla$ is the pull-back connection on $T^{\ast}M \otimes {\phi^{-1}}TX$.
\newline Using the frame adapted to the almost complex structure on $TN$, Equation~\eqref{glob} can be written:
\begin{align}
&\partial_{i{\overline{\jmath}}}\phi^{\ssc{A}} - {}^{\ssc N}\Gamma^{\ssc{C}}_{i\overline{\jmath}} \partial_{\ssc{C}}\phi^{\ssc{A}} + {}^{\ssc X}\Gamma^{\ssc{A}}_{\ssc{JL}} \partial_{i}\phi^{\ssc{J}}\partial_{\overline{\jmath}}\phi^{\ssc{L}} = \, 0 \label{IVequa21a}\\
&\forall A\in \{1,\dots ,r \} \mbox{ and } \forall i,j\in \{1,\dots ,n\}, \notag
\end{align}
where we use the notations
$$\partial_{\ssc{C}}\phi^{\ssc{A}} \quad \text{ for } \quad
d\phi^{\ssc{A}}(\vartheta_{\ssc{C}}) 
\quad \text{ and }\quad 
\partial_{i{\overline{\jmath}}}\phi^{\ssc{A}} \quad \text{ for } \quad
d(d\phi^{\ssc{A}}({\vartheta_{i}}))(\vartheta_{\overline{\jmath}}).$$

If we suppose $N$ to be a {\he} manifold with local complex coordinates $\{y^{1},y^{\overline{1}},\dots,y^{n},y^{\overline{n}}\}$ then
locally $\phi$ satisfies:
\begin{align}
&\dif{{}^{2}\phi^{\ssc{A}}}{y^{i}\partial y^{\overline{\jmath}}} - {}^{\ssc N}\Gamma^{\ssc{C}}_{i\overline{\jmath}} \dif{\phi^{\ssc{A}}}{y^{\ssc{C}}} + {}^{\ssc X}\Gamma^{\ssc{A}}_{\ssc{JL}}\dif{\phi^{\ssc{J}}}{y^{i}}\dif{\phi^{\ssc{L}}}{y^{\overline{\jmath}}} =\, 0 \tag{11G}\label{oog}\\
&\forall A\in \{1,\dots ,r \} \mbox{ and } \forall i,j\in \{1,\dots ,n\}. \notag
\end{align}
\end{de}

\begin{prop}\cite{EelLem88}
A {\oo} map from an almost {\he} {\mfd} to a {\ri} {\mfd} is {\ha}.
\end{prop}

\begin{de}\label{IVdefoog}
A smooth map from a complex {\mfd} $N^{n}$ to a {\ri} {\mfd} $X^{r}$ satisfying:
\begin{equation}
\nabla^{(0,1)}d'\phi = 0 \label{IVequaP}
\end{equation}       
is called a {\em{{\ph} map}} (\cite{Uda88A}), where
$$\left( \nabla^{(0,1)}d'\phi \right) (\overline{Z},W) = 
{}^{\phi}\nabla_{\overline{Z}} \left( d'\phi(W) \right) - 
d'\phi \left( \overline{\partial}_{\overline{Z}} W \right) \qquad \forall Z,W\in T^{(1,0)}N^{n}.$$
A map is a {\ph} map if and only if its restriction to any holomorphic curve is harmonic (see~\cite[Prop. 1.1]{OhnVal90}).
\newline We remark that:
$$\nabla^{(0,1)}d'\phi = 0 \Leftrightarrow \nabla^{(1,0)}d''\phi = 0. $$
The expression for the left-hand side of~\eqref{IVequaP} in a local frame is:
\begin{align*}
&\left( \nabla^{(0,1)}d'\phi \right) ( \vartheta_{\overline{\jmath}} , \vartheta_{i} ) = \\
&\left( \partial_{i\overline{\jmath}}\phi^{A} + {}^{\ssc X}\Gamma^{\ssc{A}}_{\ssc{JL}} \partial_{i}\phi^{\ssc{J}}\partial_{\overline\jmath}\phi^{\ssc{L}} 
\right) \varpi_{\ssc{A}}.
\end{align*}
Hence the equation for pluriharmonicity of a map is
\begin{align} 
&\partial_{i\overline{\jmath}}\phi^{A}  + {}^{\ssc X}\Gamma^{\ssc{A}}_{\ssc{JL}} \partial_{i}\phi^{\ssc{J}}\partial_{\overline\jmath}\phi^{\ssc{L}} = \, 0 \tag{PH}\label{equaph}\\
&\forall A\in \{1,\dots ,r \} \mbox{ and } \forall i,j\in \{1,\dots ,n\}.\notag
\end{align}
\end{de}

\begin{rk}
In local frames, the difference between {\ph} maps and {\oo} maps is the extra terms of the type ${}^{\ssc N}\Gamma_{i\overline\jmath}^{\ssc{C}}\partial_{\ssc{C}}\phi^{\ssc{I}}$, which appear for {\oo} maps.
\newline If we suppose $N^{n}$ to be a K\"{a}hler manifold, then by Proposition~\ref{link}.3, those terms disappear and the notions of {\oo} maps and {\ph} maps coincide.

A {\ph} map on a complex {\mfd} $N$ is harmonic with respect to any K\"{a}hler metric on $N$ and within this viewpoint can be seen as an extension of the notion of harmonic maps on surfaces (see Udagawa~\cite{Uda88A,Uda94} ). 
\end{rk}

Using Proposition~\ref{link} and the expressions of {\oo}, {\ph} and $\pm$holomorphic in adapted frames, we study when those three different concepts coincide.

\begin{de}
By {\em local} we shall mean {\em defined on an open subset}.
\end{de}

\begin{prop} \label{prop3}
Let $N^{n}$ and $X^{r}$ be {\he} {\mfd}s, then the condition that all (local) $\pm$holomorphic maps from $N^{n}$ to $X^{r}$ be {\ph} is equivalent to $X^{r}$ being {\ka}. In particular, $\pm$holomorphic {\em functions} on a {\he} {\mfd} are always {\ph}.
\end{prop}
\begin{proof}
The equation for a {\ph} map is:
\begin{align*}
&\dif{{}^{2}\phi^{\ssc{A}}}{y^{i}\partial y^{\overline{\jmath}}}  + {}^{\ssc X}\Gamma^{\ssc{A}}_{\ssc{JL}}
\dif{\phi^{\ssc{J}}}{y^{i}}
\dif{\phi^{\ssc{L}}}{y^{\overline{\jmath}}}
=0 \\
&\forall A \in \{1,\overline{1},\dots ,r,\overline{r} \} \mbox{ and } \forall i,j\in \{1,\dots ,n\}.
\end{align*} 
For a holomorphic map this reduces to the following equation:
\begin{align}
{}^{\ssc X}\Gamma^{\ssc{A}}_{m\overline{p}}
\dif{\phi^{m}}{y^{i}}
\dif{\phi^{\overline{p}}}{y^{\overline{\jmath}}} 
&=0 \qquad \forall A \in \{1,\overline{1},\dots ,r,\overline{r} \} \mbox{ and } \forall i,j\in \{1,\dots ,n\},\label{eqc1}
\end{align}

By choosing holomorphic maps with prescribed first derivatives at a point, 
Equation~\eqref{eqc1} shows that ${}^{\ssc X}\Gamma^{\ssc{A}}_{m\overline{p}} =0 \quad \forall A,m,p$, so that, by Proposition~\ref{link}, $X$ is (1,2)-symplectic and so {\ka}.
\newline The converse statement comes directly from the definitions.
\end{proof}

Now we compare {\oo} and $\pm$holomorphic.

\begin{prop}  \label{hoog}
Let $N^{n}$ and $X^{r}$ be Hermitian manifolds, then the condition that {\em all} (local) $\pm$holomorphic maps from $N$ to $X$ be {\oo} is equivalent to the condition that $N^{n}$ and $X^{r}$ be {\ka}.
\end{prop}
\begin{proof}
We equip $N^{n}$ with the local coordinates 
$(y^{1}, y^{\overline{1}},\dots, y^{n},  y^{\overline{n}})$.
\newline Let $\phi :N^{n}\rightarrow X^{r}$ be a holomorphic map, i.e.
$$ \dif{\phi^{i}}{y^{\overline{\jmath}}} = 0 \quad \forall i=1,\ldots,r, \, \forall j=1,\ldots,n.$$
We recall that $\phi$ is a {\oo} map if and only if:
\begin{align}  \label{ab}
&\dif{{}^{2}\phi^{\ssc{A}}}{y^{i}\partial y^{\overline{\jmath}}} - {}^{\ssc N}\Gamma^{\ssc{C}}_{i\overline{\jmath}} \dif{\phi^{\ssc{A}}}{y^{\ssc{C}}} + {}^{\ssc X}\Gamma^{\ssc{A}}_{\ssc{JL}}\dif{\phi^{\ssc{J}}}{y^{i}}\dif{\phi^{\ssc{L}}}{y^{\overline{\jmath}}}
 =0 \\
&\forall \, A\in \{1,\overline{1},\dots ,r,\overline{r} \} \mbox{ and } \forall \, i,j\in \{1,\dots ,n\}. \notag
\end{align}
As we supposed $\phi$ to be holomorphic, Equation~\eqref{ab} is equivalent to:
\begin{align}
- {}^{\ssc N}\Gamma^{k}_{i\overline{\jmath}} \dif{\phi^{l}}{y^{k}} + 
{}^{\ssc X}\Gamma^{l}_{m\overline{p}}
\dif{\phi^{m}}{y^{i}}  \dif{\phi^{\overline p}}{y^{\overline{\jmath}}}
 &=0 \qquad \forall \, l\in \{1,\dots ,r \} \mbox{ and } \forall \, i,j\in \{1,\dots ,n\}, \label{eqa}
\end{align}
The proof consists of three steps: first we show that ${}^{\ssc N}\Gamma^{k}_{i\overline{\jmath}} = 0 \quad \forall i,j,k$ such that $k\neq i$ or $k\neq j$, then we establish that ${}^{\ssc X}\Gamma^{l}_{m\overline p} = 0 \quad \forall l,m,p$, and finally we prove that ${}^{\ssc N}\Gamma^{k}_{i\overline\jmath} = 0 \quad \forall i,j,k$. We work at a point of $N$.
\newline Step 1) Fix the indices $l,i$ and $j$. Equation~\eqref{eqa} is:
$$- {}^{\ssc N}\Gamma^{k}_{i\overline{\jmath}} \dif{\phi^{l}}{y^{k}} + 
{}^{\ssc X}\Gamma^{l}_{m\overline{p}}
\dif{\phi^{m}}{y^{i}}  \dif{\phi^{\overline p}}{y^{\overline{\jmath}}}
 =0 ,$$
so choosing a {\ho} map $\phi$ such that $\dif{\phi^{\ssc A}}{y^{\ssc I}} = \delta^{\ssc{A}l}_{\ssc{I}k_{0}}$ with $k_{0}\neq i$ or $k_{0}\neq j$ (possible if $n > 1$), \eqref{eqa} becomes:
$${}^{\ssc N}\Gamma^{k_{0}}_{i\overline{\jmath}} = 0,$$
which proves Step 1. We remark that if $n=1$, any {\he} metric on $N$ is {\ka} (cf.~\cite{Kobayashi}) so ${}^{\ssc N}\Gamma^{k}_{i\overline{\jmath}} = 0 \quad \forall i,j,k$.
\newline Step 2) Consider Equation~\eqref{eqa} when $i\neq j$. In this case, whichever value $k$ takes, either $k\neq i$ or $k\neq j$. Because of Step 1) (and the remark for $n=1$), with those indices, \eqref{eqa} reads:
$${}^{\ssc X}\Gamma^{l}_{m\overline{p}}
\dif{\phi^{m}}{y^{i}}  \dif{\phi^{\overline p}}{y^{\overline{\jmath}}} =0.$$
Fix $m_{0}, p_{0}$ and choose $\phi$ {\ho} such that $\dif{\phi^{\ssc M}}{y^{\ssc I}} = \delta^{\ssc{M}m_{0}}_{\ssc{I}i} + \delta^{\ssc{M}p_{0}}_{\ssc{I}j}$. For this map Equation~\eqref{eqa} yields
$${}^{\ssc X}\Gamma^{l}_{m_{0}\overline{p_{0}}} = 0,$$
and this shows Step 2).
\newline Step 3) Because of Step 2) Equation~\eqref{eqa} now reads:
$${}^{\ssc N}\Gamma^{k}_{i\overline{\jmath}} \dif{\phi^{l}}{y^{k}} = 0.$$
Choosing a holomorphic map such that $\dif{\phi^{\ssc M}}{y^{\ssc I}} = \delta^{\ssc{M}l_{0}}_{\ssc{I}k_{0}}$ proves that
$${}^{\ssc N}\Gamma^{k_{0}}_{i\overline{\jmath}} = 0 \quad \forall k_{0},i,j$$
and in particular that
$${}^{\ssc N}\Gamma^{k}_{k\overline k} = 0 \quad \forall k,$$
ending Step 3).
This proves that if $\pm$holomorphic maps are {\oo} then $X^{r}$ and $N^{n}$ must be (1,2)-symplectic manifolds and therefore {\ka}. 
\newline The converse assertion is clear from the definitions of (1,2)-symplectic manifolds, {\oo} and $\pm$holomorphic maps.
\end{proof}

We now give an example where the notions of {\oo} and {\ph} are distinct.
\begin{ex}
The Hopf manifold (cf.~\cite[Vol. II]{Kobayashi}) is defined as $H = \left( {\cn}^{\,n} - \{ 0 \} / \Delta_{\lambda} \right)$ $(n>1)$, where $\Delta_{\lambda}$ is the group generated by the transformation $ z\mapsto \lambda z , \quad z\in {\cn} - \{0,1\}$.
\newline $H$ is diffeomorphic to ${\sn}^{1}\times{\sn}^{2n-1}$, and is a compact {H}ermitian manifold equipped with the metric
$$ds^{2} = \frac{1}{ \sum^{n}_{k=1} z^{k}z^{\overline k}} \, \sum^{n}_{j=1} dz^{j}\otimes dz^{\overline\jmath}.$$
Its fundamental K{\" a}hler form $\Omega$ is
$$\Omega = \frac{-i}{\pi\sum^{n}_{k=1} z^{k}z^{\overline k}} \, \sum^{n}_{j=1} dz^{j}\wedge dz^{\overline\jmath}.$$
The manifold $H$ is not K{\" a}hler since its second Betti number $b_{2}(H)$ is equal to $b_{2}({\sn}^{1}\times{\sn}^{2n-1})$ which is 1, whereas the Betti numbers $b_{2n}$ of a K{\" a}hler manifold are even.
\newline A simple computation shows that
\begin{align*}
\Gamma^{k}_{i\overline\jmath} =& \frac{1}{2 \mid z\mid^{2}}\left( \delta_{ij} z^{k} - \delta_{ik} z^{j} \right) \\
\Gamma^{\overline k}_{i\overline\jmath} =& \frac{1}{2 \mid z\mid^{2}}\left( \delta_{ij} z^{\overline k} - \delta_{jk} z^{\overline\imath} \right)
\end{align*}
therefore
\begin{align}
&\frac{\partial^{2} \phi}{\partial z^{i}\partial  z^{\overline\jmath}} - \Gamma^{k}_{i\overline\jmath} \frac{\partial\phi}{\partial z^{k}} - \Gamma^{\overline k}_{i\overline\jmath} \frac{\partial\phi}{\partial  z^{\overline k}} = \notag\\
&\frac{\partial^{2} \phi}{\partial z^{i}\partial  z^{\overline\jmath}} + \frac{1}{2 \mid z\mid^{2}} 
\left[ z^{j} \frac{\partial\phi}{\partial z^{i}} + z^{\overline\imath} \frac{\partial\phi}{\partial  z^{\overline\jmath}} - \delta_{ij} \sum_{k} \left( z^{\overline k} \frac{\partial\phi}{\partial  z^{\overline k}} + z^{k} \frac{\partial\phi}{\partial z^{k}} \right) \right]. \label{last}
\end{align}
If $\phi$ is a {\oo} function on $H$ then
\begin{equation} 
\frac{\partial^{2} \phi}{\partial z^{i}\partial  z^{\overline\jmath}} - \Gamma^{k}_{i\overline\jmath} \frac{\partial\phi}{\partial z^{k}} - \Gamma^{\overline k}_{i\overline\jmath} \frac{\partial\phi}{\partial z^{\overline k}} = \quad 0 \qquad \forall i,j = 1,\ldots,n. \tag{**}\label{star2}
\end{equation}
In the case where $\phi$ is holomorphic Equation~\eqref{star2}, for $i\neq j$, becomes:
$$  z^{j} \frac{\partial\phi}{\partial z^{i}} = 0,$$
showing that $\phi$ has to be constant.
\end{ex}

We have shown
\begin{prop}
Any local $\pm${\ho} {\oo} function on $H$ is constant.
\end{prop}

\begin{rk} \label{propIV5}
Since $H$ is {\he}, any (non-constant) $\pm${\ho} function on $H$ is {\ph}, therefore local {\ph} functions exist.
\end{rk}

Finally we compare {\oo} maps and {\ph} maps on Hermitian manifolds and the situation is what the last two propositions suggest:

\begin{prop} \label{propIV6}
Let $N^{n}$ be a Hermitian manifold and $X^{r}$ a {\ka} {\mfd}, then the condition ``a map from $N^{n}$ to $X^{r}$, is a (local) {\oo} map {\iff} it is {\ph}'' is equivalent to $N$ being {\ka}.
\end{prop}
\begin{proof}
The definitions of {\oo} map and {\ph} map for $\phi:N^{n}\rightarrow X^{r}$ coincide when:
\begin{equation}
{}^{\ssc N}\Gamma^{\ssc C}_{i\overline{\jmath}} 
\dif{\phi^{\ssc{I}}}{y^{\ssc C}} = 0 \qquad \forall I\in \{1,\overline{1},\dots ,r,\overline{r} \} \quad \forall i,j\in \{1,\dots ,n\}. \label{IVstar}
\end{equation}
As $X^{r}$ is {\ka}, $\pm${\ho} maps from $N^{n}$ to $X^{r}$ are {\ph} (Proposition~\ref{prop3}). We choose $\phi$ {\ho} and such that $$\frac{\partial\phi^{\ssc A}}{\partial y^{\ssc C}} = 
\delta^{{\ssc A}\alpha}_{{\ssc C}k_{0}},$$ 
then Equation~\eqref{IVstar} reads:
$$ {}^{\ssc N}\Gamma^{k_{0}}_{i\overline{\jmath}} =0 \quad \forall i,j \in \{1,\dots,n\}.$$
Since this holds for any $k_{0} = 1,\ldots,n$, $N$ is (1,2)-symplectic and so {\ka} since integrable.
\end{proof}

\begin{rk}
In Propositions~\ref{prop3},~\ref{hoog} and~\ref{propIV6} the need for the integrability of the complex structures on $N^{n}$ and $X^{r}$ is in the use of $n$ linearly independent holomorphic functions as ``test functions'', the existence of these being equivalent to the integrability of the structures.
\end{rk}

\begin{prop}\cite{OhnVal90}
Let $(N,J,g)$ be a {\he} {\mfd} and $(X,h)$ an almost {\he} {\mfd}. 
Then $(N,J,g)$ is cosymplectic {\iff} any {\ph} map $\phi : (N,J,g) \to (X,h)$ is {\ha}.
\end{prop}

\begin{de}
We shall call $\psi$ a {\em {\ph} morphism} if it pulls back (local) {\ph} functions to (local) {\ph} functions, i.e. for all functions $\phi$ satisfying 
$\nabla^{(0,1)}d' \phi = 0$ we have
$$\nabla^{(0,1)}d'(\phi\circ\psi ) = 0 .$$
\end{de}

\begin{rk}
Since a complex-valued function is {\ph} {\iff} its real and imaginary parts are real-valued {\ph} functions, we can choose to consider, in the above definition, real-valued or complex-valued functions.
\end{rk}

\begin{lem} \label{lem1}
Let $N$ be a {\he} {\mfd}.
Given $p\in N$ and any constants $(C_{{\ssc A}})_{{\ssc{A=1,\overline{1},\dots ,n,\overline{n}}}}$ and 
$(C_{{\ssc{AB}}})_{{\ssc{A,B=1,\overline{1},\dots ,n,\overline{n},typeA\neq typeB}}}$ such that $C_{{\ssc{AB}}}=C_{{\ssc{BA}}}$, there exists a (complex-valued) {\ph} function, $h$, defined on a neighbourhood of $p$ satisfying 
$\dif{h}{y^{\ssc{A}}}(p) = C_{\ssc{A}}$
and
$\dif{{}^{2}h}{y^{\ssc{A}}\partial y^\ssc{B}}(p) = C_{\ssc{AB}}.$
\end{lem}

\begin{proof}
We choose complex local coordinates $(y^{1}, y^{\overline{1}},\dots, y^{n}, y^{\overline{n}})$ in a neighbourhood of the point $p \in N$.
\newline We shall look for $\pm${\ho} functions; in this case, without loss of generality we can take $N={\cn}^{\,n}$ and $X={\cn}$. This is then easy, for example:
$$h(z) = \sum_{{\ssc{A}} =1,\overline{1}}^{n,\overline{n}} C_{{\ssc{A}}}   y^{{\ssc{A}}} + \frac{1}{2}\sum_{{\ssc{A}},{\ssc{B}} =1,\overline{1}}^{n,\overline{n}} C_{{\ssc{A}}{\ssc{B}}}y^{{\ssc{A}}} y^{{\ssc{B}}},$$
where the $C_{{\ssc{A}}}$ and $C_{{\ssc{A}}{\ssc{B}}}$ are complex constants, is {\ho} and has the prescribed first and second derivatives.
\end{proof}

\begin{prop} \label{Vprop6}
Let $\psi$ be a smooth map from a {\he} {\mfd} $M$ to a {\he} {\mfd} $N$. Then $\psi$ is a {\ph} morphism if and only if it is a $\pm$holomorphic  {\ph} map.
\end{prop}

\begin{proof}
Let $\psi: M \to N$ be a {\ph} morphism, then it must satisfy
$$\nabla^{(0,1)} d'(\phi\circ\psi )(\overline{Z},W) = 0 $$
for all $Z,W\in T^{(1,0)}M $ and all local {\ph} functions $\phi : N\to {\cn}$.
\newline We first compute the  chain rule for the operator $\nabla^{(0,1)} d'$.
\newline Recall that:
\begin{align*}
&d'(\phi\circ\psi) = d(\phi\circ\psi)|_{T^{(1,0)}M}
= d\phi\circ d\psi |_{T^{(1,0)}M} = d\phi\circ d'\psi .
\end{align*}
Hence
\begin{align*}
&\nabla^{(0,1)} d'(\phi\circ\psi )(\overline{Z},W) = \\
&{}^{\phi\circ\psi}\nabla_{\overline{Z}} 
\left( d'(\phi\circ\psi )(W) \right) - 
d'(\phi\circ\psi ) \left( \overline{\partial}_{\overline{Z}} W \right) =\\
&{}^{\phi\circ\psi}\nabla_{\overline{Z}} 
\left( d\phi\circ d'\psi (W) \right) - 
d\phi\circ d'\psi  \left( \overline{\partial}_{\overline{Z}} W \right) =\\
&\left( {}^{\phi}\nabla_{d''\psi (\overline{Z})} 
 d\phi \right) \left(d'\psi(W) \right) + 
d\phi \left( {}^{\psi}\nabla_{\overline{Z}} \left(d\psi(W)\right)  \right) - 
d\phi\circ d'\psi  \left( \overline{\partial}_{\overline{Z}} W \right)= \\
&\left(\nabla d\phi\right) \left( d''\psi (\overline{Z}), d'\psi(W) \right) + 
d\phi \left[ {}^{\psi}\nabla_{\overline{Z}}\left( d\psi(W)\right) - 
 d'\psi  \left( \overline{\partial}_{\overline{Z}} W \right) \right]= \\
&\left[ \nabla d\phi \left( d''\psi , d'\psi \right) + 
d\phi \left( \nabla^{(0,1)} d'\psi \right) \right] (\overline{Z},W).
\end{align*}
If we choose $Z,W$ to be the vectors ${\frac{\partial}{\partial z^{i}}},{\frac{\partial}{\partial z^{j}}}$ of the canonical frame of $T^{(1,0)}M $:
\begin{align}
&\nabla^{(0,1)} d'(\phi\circ\psi ) \left({\frac{\partial}{\partial z^{\overline \jmath}}} , {\frac{\partial}{\partial z^{i}}}\right)= \notag \\
&\dif{\phi}{y^{\ssc{C}}}
\left[ \dif{{}^{2}\psi^{\ssc{C}}}{z^i\partial z^{\overline{\jmath}}} 
 + {}^{\ssc N}\Gamma^{\ssc{C}}_{\ssc{JL}}
\dif{\psi^{\ssc{J}}}{z^{i}}\dif{\psi^{\ssc{L}}}{z^{\overline{\jmath}}}
\right] +  \label{equa}\\ 
&
\sum_{\substack{\ssc{A,B}=1,\overline{1}  \\\ssc{A},\ssc{B}\text{ of same type}}}
^{n,\overline{n}}
\dif{\psi^{\ssc{A}}}{z^{i}}
\dif{\psi^{\ssc{B}}}{z^{\overline{\jmath}}}
\left[
\dif{{}^{2}\phi}{y^{\ssc{A}}\partial y^{\ssc{B}}}  
 - {}^{\ssc N}\Gamma^{\ssc{C}}_{\ssc{AB}} \dif{\phi}{y^{\ssc{C}}} \right]. \notag 
\end{align}
The ``only if'' part of the proposition follows from~\eqref{equa}.
\newline Conversely, by Lemma~\ref{lem1} we can choose various different values for $\dif{\phi}{y^{\ssc{C}}}$ and 
$\dif{{}^{2}\phi}{y^{\ssc{B}}\partial y^{\ssc{A}}}$ at a given point:
\newline If we take $\dif{{}^{2}\phi}{y^{\ssc{B}}\partial y^{\ssc{A}}}= 0$ for all $A,B$ and
$ \dif{\phi}{y^{\ssc{C}}}= 0$ for all $C$ except for $l_{0}$ (we will suppose $ \dif{\phi}{y^{\ssc{l_{0}}}}= 1$) then Equation~\eqref{equa} becomes:
\begin{align}
&\dif{{}^{2}\psi^{l_{0}}}{z^{i}\partial z^{\overline{\jmath}}} +  \sum_{\ssc{J,L}=1,\overline{1}}^{n,\overline{n}}
{}^{\ssc N}\Gamma^{l_{0}}_{\ssc{JL}}
\dif{\psi^{\ssc{J}}}{z^{i}}
\dif{\psi^{\ssc{L}}}{z^{\overline{\jmath}}} 
- \sum_{\substack{\ssc{A,B}=1,\overline{1}  \\\ssc{A},\ssc{B}\text{ of same  type}}}
^{n,\overline{n}}
{}^{\ssc N}\Gamma^{l_{0}}_{\ssc{AB}}
\dif{\psi^{\ssc{A}}}{z^{i}}
\dif{\psi^{\ssc{B}}}{z^{\overline{\jmath}}}
= \notag \\
&\dif{{}^{2}\psi^{l_{0}}}{z^{i}\partial z^{\overline{\jmath}}} 
+  \sum_{\substack{\ssc{I,J}=1,\overline{1}  \\\ssc{I},\ssc{J}\text{ of different type}}}
^{n,\overline{n}}
{}^{\ssc N}\Gamma^{l_{0}}_{\ssc{JL}}
\dif{\psi^{\ssc{J}}}{z^{i}}
\dif{\psi^{\ssc{L}}}{z^{\overline{\jmath}}}
  = 0, \label{equa15}
\end{align}
for all $l_{0}\in \{1,\overline{1},\dots ,n,\overline{n}\}$.
If we take $\phi$ such that $\dif{\phi}{y^{\ssc{C}}}=0 $ for all $C\in\{1,\overline{1},\dots ,n,\overline{n}\} $ and $\dif{{}^{2}\phi}{y^{\ssc{B}}\partial y^{\ssc{A}}} = 0$ for all $A,B$ but for $A_{0},B_{0}$ with type $A_{0}$ different from type $B_{0}$ (we will suppose that $\dif{{}^{2}\phi}{y^{\ssc{B_{0}}}\partial y^{\ssc{A_{0}}}} = 1$), Equation~\eqref{equa} becomes:
\begin{equation}
\label{holo}
\dif{\psi^{\beta_{0}}}{z^{\overline{\jmath}}}\dif{\psi^{\alpha_{0}}}{z^{i}} =
 0,
\end{equation}
which we obtain for all $\alpha_{0},\beta_{0}=1,\dots ,n, \, i,j= 1,\dots,m$.
\newline Either (i) for all $\alpha=1,\dots ,n,$ and all $i$: $$\dif{\psi^{\alpha}}{z^{i}}=0 $$
\newline or (ii) there exists $\alpha_{0}$ and $i_{0}$ such that: 
$$\dif{\psi^{\alpha_{0}}}{z^{i_{0}}}\neq 0.$$
Then Equation~\eqref{holo} shows that:
$$\dif{\psi^{\alpha}}{z^{\overline{\imath}}} = 0,$$
for all $\alpha=1,\dots ,n$ and all $i$.     
\newline Therefore at each point $p \in M$, $\psi$ is either antiholomorphic or  holomorphic.
\newline We shall use an argument of Siu~\cite{Siu80} to prove that $\psi$ is either holomorphic or antiholomorphic on the whole of $M$.
\newline Let
\begin{align*}
U = \left\{ p\in M \quad \vert \quad\partial\psi (p) \equiv 0 \right\} \\
V = \left\{ p\in M \quad \vert \quad\overline{\partial}\psi (p) \equiv 0 \right\} .
\end{align*}
As 
$$ U \cup V = M,$$
either $\overset{\circ}{U}\neq\emptyset$ or $\overset{\circ}{V}\neq\emptyset$.
Suppose that $\overset{\circ}{V}\neq\emptyset$ and let $\widetilde{V}$ be the largest connected open subset of $M$ included in $V$.
\newline If $M=\widetilde{V}$ then $\psi$ is {\ho}.
\newline Assume that $M\neq \widetilde{V}$.
\newline Let $p$ be a boundary point of $V$. Let $W$ be an open connected neighbourhood of $p$ in $M$ such that:
\begin{enumerate}
\item there exists a {\ho} coordinate system $(z^{i})$ on some open neighbourhood of $\overline{W}$,
\item there exists a {\ho} coordinate system $(w^{\alpha})$ on some open neighbourhood of $\overline{\psi(W)}$.
\end{enumerate}
Since $\psi$ satisfies Equation~\eqref{equa15}:
\begin{align} \label{ho2}
&\dif{{}^{2}\psi^{\alpha}}{z^{i}\partial z^{\overline{\jmath}}} +  \sum_{\substack{\ssc{I,J}=1,\overline{1}  \\\ssc{I},\ssc{J}\text{ of different type}}}
^{n,\overline{n}}
{}^{\ssc N}\Gamma^{\alpha}_{\ssc{JL}}
\dif{\psi^{\ssc{J}}}{z^{i}}
\dif{\psi^{\ssc{L}}}{z^{\overline{\jmath}}} =0 \\
&\qquad \forall \alpha\in \{1,\dots ,r \} \mbox{ and } \forall i,j\in \{1,\dots ,n\}. \notag
\end{align} 
Since any partial derivative of any term of Equation~\eqref{ho2} is bounded on $W$, applying the operator $\partial_{\overline{k}}$ to the trace of \eqref{ho2} yields:
\begin{equation}
\left| {}^{\ssc M}\Delta \left(\dif{\psi^{\alpha}}{z^{\overline{k}}} \right) \right| \leq
C \left( \sum_{i,j} 
\left| \dif{{}^{2}\psi^{\alpha}}{z^{i}\partial z^{\overline{\jmath}}} \right| + 
 \sum_{j,L} \left| \dif{\psi^{\ssc{L}}}{z^{\overline{\jmath}}} \right| +
\sum_{j,J} \left| \dif{{}^{2}\psi^{\ssc{J}}}{z^{\overline k}\partial z^{\overline{\jmath}}} \right| \right),
\end{equation}
where ${}^{\ssc M}\Delta = -2 \sum_{i,j} g^{i\overline{\jmath}} 
\dif{{}^{2}}{z^{i}\partial z^{\overline{\jmath}}}$
and $C$ is a positive constant.
\newline Set 
$$ \dif{\psi^{\alpha}}{z^{\overline{k}}} = u^{\alpha}_{\overline{k}} + i v^{\alpha}_{\overline{k}}.$$ 
Then Equation~\eqref{ho2} means that:
$$\left| {}^{\ssc M}\Delta u^{\alpha}_{\overline{k}} \right|^{2}  + 
\left| {}^{\ssc M}\Delta v^{\alpha}_{\overline{k}} \right|^{2} \leq C' \sum_{\beta,j} \left\{ 
\left| \g( u^{\beta}_{\overline{\jmath}} ) \right|^{2}  +
\left| \g( v^{\beta}_{\overline{\jmath}} ) \right|^{2}  +
\left|  u^{\beta}_{\overline{\jmath}}  \right|^{2}  +
\left|  v^{\beta}_{\overline{\jmath}}  \right|^{2} \right\}
$$
On $W$. Applying Aronszajn's unique continuation Theorem~\cite{Aro57} shows that since $u^{\beta}_{\overline{\jmath}} = v^{\beta}_{\overline{\jmath}} = 0$ on $W \cap \widetilde{V}$ then $u^{\beta}_{\overline{\jmath}} = v^{\beta}_{\overline{\jmath}} = 0$ on $\widetilde{V}$. But $\widetilde{V}$ was supposed to be maximal therefore $\psi$ is {\ho}.
\newline The anti-{\ho} case is similar.
\newline Now note that for a $\pm${\ho} map, Equation~\eqref{equa15} becomes the equation of {\ph}ity.
\end{proof}

\begin{prop} \label{Vprop7}
The set of pluriharmonic morphisms from a {\he} {\mfd} to a {\ka} {\mfd} is exactly the set of $\pm${\ho} maps.
\end{prop}
\begin{proof}
We have shown that the {\ph} morphisms between {\he} {\mfd}s are precisely the {\ph} maps which are $\pm${\ho}, but because of Proposition~\ref{prop3}, when the target is {\ka} these are exactly the $\pm${\ho} maps.
\end{proof}

\begin{prop}
A {\ph} morphism from a {\he} {\mfd} to a {\he} {\mfd} pulls back (local) {\ph} {\em maps} to (local) {\ph} {\em maps}.
\end{prop}
\begin{proof}
Let $\psi : M\rightarrow N$ be a {\ph} morphism (i.e. a $\pm${\ho} {\ph} map) and $\phi : N\rightarrow P$ a {\ph} map.
\newline The map $\phi\circ\psi$ is {\ph} if 
$$\nabla^{(0,1)} d'(\phi\circ\psi )(\overline{Z},W) = 0 \qquad \forall Z,W\in T^{(1,0)}M. $$
But
\begin{align}
&\nabla^{(0,1)} d'(\phi\circ\psi )(\overline{Z},W) =\notag\\ 
&d\phi\circ\nabla^{(0,1)} d'\psi (\overline{Z},W) + \nabla d\phi (d\psi (\overline{Z}), d\psi (W) )= \quad \text{ as shown in Proposition~\ref{Vprop6} } \notag\\
&\nabla d \phi (d\psi (\overline{Z}), d\psi (W) ) =\notag\\
&\nabla_{d\psi (\overline{Z})} d\phi (d\psi (W)) - d\phi \left( \nabla_{d\psi (\overline{Z})} d\psi (W) \right). \label{equat}
\end{align}
Suppose that $\psi$ is holomorphic then $d\psi (W)\in T^{(1,0)}N$, therefore 
$$d\phi (d\psi (W)) = d'\phi (d\psi (W))$$ 
and, as the bundle $T^{(1,0)}N \rightarrow N$ is {\ho} ,
$$ \nabla_{d\psi (\overline{Z})} d\psi (W)  =   \overline{\partial}_{d\psi (\overline{Z})} d\psi (W).$$
Therefore the right-hand side of \eqref{equat} can be written
$$\nabla^{(0,1)} d'\phi(d\psi (\overline{Z}),d\psi (W)),$$
which vanishes since $\phi$ is {\ph}.
\end{proof}


\begin{thebibliography}{BFPP93}

\bibitem[Aro57]{Aro57}
N.~Aronszajn.
\newblock {A} unique continuation theorem for solutions of elliptic partial
  differential equations or inequalities of second order.
\newblock {\em J. Math.\ Pures Appl.}, 9:235--249, 1957.

\bibitem[BFPP93]{BurFerPed93}
F.E. Burstall, D.~Ferus, F.~Pedit, and U.~Pinkall.
\newblock {H}armonic tori in symmetric spaces and commuting {H}amiltonian
  systems on loop algebras.
\newblock {\em Ann.\ of Math.}, 138:173--212, 1993.

\bibitem[EL78]{EelLem78}
J.~Eells and L.~Lemaire.
\newblock {A} report on harmonic maps.
\newblock {\em Bull.\ London Math.\ Soc.}, 10:1--68, 1978.

\bibitem[EL88]{EelLem88}
J.~Eells and L.~Lemaire.
\newblock {A}nother report on harmonic maps.
\newblock {\em Bull.\ London Math.\ Soc.}, 20:385--524, 1988.

\bibitem[Fug78]{Fug78}
B.~Fuglede.
\newblock {H}armonic morphisms between {R}iemannian manifolds.
\newblock {\em Ann.\ Inst.\ Fourier (Grenoble)}, 28:107--144, 1978.

\bibitem[GH80]{GrayHerv}
A.~Gray and L.~M. Hervella.
\newblock The {S}ixteen {C}lasses of {A}lmost {H}ermitian {M}anifolds and their
  {L}inear {I}nvariants.
\newblock {\em Ann. Mat. Pura Appl.}, 123:35--58, 1980.

\bibitem[Ish79]{Ish79A}
T.~Ishihara.
\newblock {A} mapping of {R}iemannian manifolds which preserves harmonic
  functions.
\newblock {\em J. Math.\ Kyoto Univ.}, 19:215--229, 1979.

\bibitem[KN69]{Kobayashi}
S.~Kobayashi and K.~Nomizu.
\newblock {\em Foundations of Differential Geometry, I and II}, volume~15 of
  {\em Interscience Tracts in Pure and Applied Mathematics}.
\newblock Interscience, 1963,1969.

\bibitem[Lic70]{Lic70}
A.~Lichnerowicz.
\newblock {A}pplications harmoniques et vari\'et\'es {K}\"ahl\'eriennes.
\newblock {\em Sympos.\ Math.}, pages 341--402, 1970.

\bibitem[OV90]{OhnVal90}
Y.~Ohnita and G.~Valli.
\newblock {P}luriharmonic maps into compact {L}ie groups and factorization into
  unitons.
\newblock {\em Proc.\ London Math.\ Soc.}, 61:546--570, 1990.

\bibitem[Sal85]{Sal85}
S.M. Salamon.
\newblock {H}armonic and holomorphic maps.
\newblock In E.~Vesentini, editor, {\em {G}eometry {S}eminar ``{L}uigi
  {B}ianchi'' {I}{I}---1984}, volume 1164 of {\em Lecture Notes in Math.},
  pages 161--224. Springer, Berlin, Heidelberg, New York, 1985.

\bibitem[Siu80]{Siu80}
Y.T. Siu.
\newblock {T}he complex-analyticity of harmonic maps and the strong rigidity of
  compact {K}\"ahler manifolds.
\newblock {\em Ann.\ of Math.}, 112:73--111, 1980.

\bibitem[Uda88]{Uda88A}
S.~Udagawa.
\newblock {P}luriharmonic maps and minimal immersions of {K}{\"a}hler
  manifolds.
\newblock {\em J. London Math.\ Soc.}, 37:375--384, 1988.

\bibitem[Uda94]{Uda94}
S.~Udagawa.
\newblock {C}lassification of pluriharmonic maps from compact complex manifolds
  with positive first {C}hern class into complex {G}rassmann manifolds.
\newblock {\em T\^ohoku Math.\ J.}, 46:367--391, 1994.

\end{thebibliography}
\end{document}